\documentclass[12pt]{amsart}

\usepackage{psfrag}
\usepackage{graphicx}
\usepackage{pinlabel}
\usepackage{graphicx}
\usepackage{enumerate}
\usepackage{amsmath}
\usepackage{amssymb}
\usepackage{amscd}
\usepackage{picinpar}
\usepackage{times}
\usepackage{pb-diagram}
\usepackage{graphicx}
\usepackage{wrapfig}
\usepackage{xspace}
\usepackage{hyperref}
\usepackage[dvipsnames]{xcolor}
\hypersetup{hidelinks}
\hypersetup{
	colorlinks,
	citecolor=blue,
	linkcolor=red,
	urlcolor=blue}
\usepackage{wasysym}
\usepackage{overpic}
\usepackage{enumitem}
\setlist{itemsep=3pt plus 1pt minus 1pt}

\usepackage{url}
\newtheorem{theorem}{Theorem}[section]
\newtheorem{lemma}[theorem]{Lemma}
\newtheorem{proposition}[theorem]{Proposition}

\newtheorem*{theorem*}{Theorem}
\newcommand{\dmo}{\DeclareMathOperator}
\dmo{\Id}{Id}
\dmo{\id}{id}
\dmo{\pt}{pt}
\dmo{\ab}{ab}
\dmo{\GL}{GL}
\dmo{\SO}{SO}
\dmo{\SU}{SU}
\dmo{\SL}{SL}
\dmo{\PSL}{PSL}
\dmo{\im}{Im}
\dmo{\re}{Re}
\dmo{\Ad}{Ad}
\dmo{\ad}{ad}
\dmo{\tr}{tr}
\dmo{\Log}{Log}
\dmo{\Homeo}{Homeo}
\dmo{\Vol}{Vol}
\dmo{\supp}{supp}
\dmo{\Isom}{Isom}
\dmo{\Hom}{Hom}
\dmo{\Ext}{Ext}
\dmo{\Tor}{Tor}
\dmo{\PP}{\mathbf{P}}
\dmo{\Area}{Area}
\dmo{\Length}{Length}
\dmo{\inte}{int}
\dmo{\sstar}{star}
\dmo{\diam}{diam}

\newcommand{\R}{\mathbb R}

\renewcommand{\tr}{\mathrm {tr}}

\newcommand{\from}{\colon}

\usepackage{mathtools}
\mathtoolsset{showonlyrefs}

\begin{document}

\title{Spaces of Geodesic Triangulations Are Cells}

\author{Yanwen Luo}

\address{School of Mathematics and Physics, University of Science and Technology Beijing, Beijing, 100083}
\email{yanwen.luo@ustb.edu.cn}

\author{Zhenghao Rao}
\address{Department of Mathematics, Rutgers University, 110 Frelinghuysen Road, Piscataway, NJ, 08854}
\email{zhenghao.rao@rutgers.edu}

\author{Tianqi Wu}

% \address{Clark University, 950 Main St, Worcester, MA 01610, 
% }

\email{mike890505@gmail.com}

\author{Xiaoping Zhu}

\address{Department of Mathematics, Rutgers University, Piscataway, NJ 08854}

\email{xz349@math.rutgers.edu}



\begin{abstract}
It has been shown that spaces of geodesic triangulations of closed negatively curved surfaces are contractible. Here we prove that these spaces are homeomorphic to Euclidean spaces $\mathbb{R}^n$.

\end{abstract}
\maketitle

\section{Introduction}
Geodesic triangulations provide natural discrete representations of Riemannian surfaces. Once the combinatorial triangulation and its homotopy class are fixed, all such representations on a surface   form a finite-dimensional manifold, called a \textit{space of geodesic triangulations} of the surface.

These spaces are closely related to the classical groups of surface diffeomorphisms.  Smale first established the contractibility theorem for the diffeomorphism group of the disk relative to its boundary~\cite{smale1959diffeomorphisms}.  Earle--Eells \cite{earle1969fibre} later showed that 
the identity component of 
the group of surface diffeomorphisms on a surface of constant curvature is homotopy equivalent to the subgroup of orientation-preserving isometries isotopic to the identity. In particular, this implies that the identity component of the group of surface diffeomorphisms is homotopy equivalent to $\SO(3)$ for the unit 2-sphere, to the torus of translations for a flat torus, and to a point for a closed hyperbolic surface.

Spaces of geodesic triangulations of surfaces can be regarded as finite-dimensional analogues of groups of surface diffeomorphisms. 
Intuitively, if the triangulations are chosen to be sufficiently dense, we can approximate a large class of diffeomorphisms with sufficient accuracy. 
Therefore,  it is natural to ask whether a space of geodesic triangulations of a surface is homotopy equivalent to the group of diffeomorphisms of the surface. This analogy was formulated as a conjecture by Connelly--Henderson--Ho--Starbird~\cite{connelly1983problems} for a closed orientable surface equipped with a constant-curvature metric: \textit{each space of geodesic triangulations is homotopy equivalent to the group of diffeomorphisms of the surface}. Specifically, the conjecture predicts that spaces of geodesic triangulations of the unit 2-sphere, flat tori, and hyperbolic surfaces are homotopy equivalent to $\SO(3)$, a torus, and a point, respectively.

Beginning with the work of Cairns~\cite{cairns1944isotopic}, the results for convex polygons provides the first evidence for this conjecture.  Bloch--Connelly--Henderson~\cite{bloch1984space} proved that the space of geodesic triangulations of a convex polygon is homeomorphic to a Euclidean space, as a discrete counterpart of Smale's theorem. In fact, Bloch--Connelly--Henderson observed that Smale's theorem can be recovered from their discrete result by a limiting argument.
Substantial progress on the Connelly--Henderson--Ho--Starbird conjecture has been made in recent years. For flat tori, the conjecture was proved independently using two different approaches~\cite{erickson2021planar,luo2021deformation2}. For closed surfaces of negative curvature, including hyperbolic surfaces as special cases, each space of geodesic triangulations  was proved to be contractible~\cite{luo2021deformation}. On the unit 2-sphere, the space of Delaunay triangulations, if nonempty, is homotopy equivalent to $\SO(3)$~\cite{luo2022deformation}, but the homotopy type of the full space of geodesic triangulations remains unknown. Thus, the spherical case is the remaining open case of the conjecture, with only partial results in \cite{awartani1987spaces}.

Spaces of geodesic triangulations are also directly connected with the graph morphing problem in computational geometry. A morph between two geodesic triangulations of the same graph on a surface is a path in the space of geodesic triangulations.  Consequently, connectivity of these spaces guarantees the existence of a morph, while a deeper understanding of their topology and geometry provides a foundation for effective morphing algorithms. For planar graphs in the Euclidean plane, various effective algorithms have been developed~\cite{planar, floater1999morph}.  Related constructions have been generalized to flat tori~\cite{chambers2021morph, erickson2021planar}. The generalized Tutte embedding theorem for negatively curved surfaces~\cite{luo2021deformation} supplies the analogous parameterization on higher-genus hyperbolic surfaces, leading to a morphing algorithm in \cite{luo2026morphing}. In contrast to the planar, toroidal, and hyperbolic cases above, developing an efficient algorithm for morphing arbitrary geodesic triangulations on the \(2\)-sphere remains an open problem~\cite{erickson2025shelling}.

The Connelly--Henderson--Ho--Starbird conjecture concerns the homotopy type of spaces of geodesic triangulations. It does not, however, determine the topology of these spaces at infinity: an open contractible manifold need not be homeomorphic to Euclidean space, as illustrated by the complicated topology at infinity of the famous Whitehead manifold.

The current work, as a continuation of \cite{luo2021deformation}, rules out the potential pathological behavior at infinity and establishes the analogy with the Bloch--Connelly--Henderson theorem. We prove

\begin{theorem}
\label{main}
Every space of geodesic triangulations of a closed negatively curved surface is homeomorphic to some Euclidean space.
\end{theorem}
\subsection{Setup of Theorem \ref{main}} We first recall the following notation from \cite{luo2021deformation}. 
Let $M$ be a connected closed orientable smooth surface with a smooth Riemannian metric $g$ of negative Gaussian curvature.  A triangulation of $M$ is a homeomorphism $\psi$ from $|T|$ to $M$, where $|T|$ is the carrier of a 2-dimensional simplicial complex $T=(V,E,F)$ with the vertex set $V$, the edge set $E$, and the face set $F$. 
For a vertex $v\in V$, let 
% $\sstar(v):=\{\sigma\in F:v\in \sigma\}$ 
$\sstar(v):=\cup_{v\in\sigma\in F} \sigma$ 
be the closed $1$-ring neighborhood of $v$ in $T$.
We label vertices of $T$ as $1,2,...,n$, where $n=|V|$ is the number of vertices.
The edge in $E$ determined by vertices $i$ and $j$ is denoted by $ij$. Each edge is identified with the closed interval $[0,1]$. We use $N_v$ to denote the neighbors of a vertex $v$.

Let $T^{(1)}$ be the 1-skeleton of $T$. A \textit{geodesic triangulation} in the homotopy class of $\psi$ is an embedding $\varphi\from T^{(1)}\rightarrow M$ satisfying the following conditions:
\begin{enumerate}
	\item The restriction $\varphi_{ij}$ of $\varphi$ to the edge $ij$ is a geodesic parameterized at constant speed;
	\item $\varphi$ is homotopic to $\psi|_{T^{(1)}}$.
\end{enumerate}
The space of all geodesic triangulations homotopic to $\psi$ is called the \emph{space of geodesic triangulations homotopic to $\psi$}, denoted by $X(M,T,\psi)$. 
In the rest of the paper, we fix a triangulation $T$ and a homotopy class determined by $\psi$ unless stated otherwise.
For simplicity, we let $X=X(M,T,\psi)$ and call it a \emph{space of geodesic triangulations}.
% Denote by $X=X(M,T,\psi)$ the space of geodesic triangulations homotopic to $\psi$. 
Then $X$ is a metric space, with the distance function
\begin{equation}\label{eq:distance-func}
    d_X(\varphi,\phi)=\max_{x}d_{g}(\varphi(x),\phi(x)).
\end{equation}
\begin{theorem}[{\cite[Theorem 1.1]{luo2021deformation}}]
\label{contractible}
The space $X(M,T,\psi)$ is contractible. 
\end{theorem}

We restate Theorem \ref{main} more specifically as follows. 
\begin{theorem}
\label{main2}
The space $X(M,T,\psi)$ is homeomorphic to $\mathbb{R}^{2n}$, where $n = |V|$. 
\end{theorem}

We use the following notation in the rest of this paper. 
Let $\tilde X=\tilde X (M,T,\psi)$ be the ambient space containing $X$, consisting of all continuous maps $\varphi\from T^{(1)}\rightarrow M$ satisfying the following conditions:
\begin{enumerate}
	\item The restriction $\varphi_{ij}$ of $\varphi$ to the edge $ij$ is a geodesic parameterized at constant speed;
	\item $\varphi$ is homotopic to $\psi|_{T^{(1)}}$.
\end{enumerate}
We note that the elements of $\tilde X$, which are called \textit{geodesic mappings}, need not be embeddings from $T^{(1)}$ to $M$. The space $\tilde{X}$ is also a metric space, equipped with the same distance function given by \eqref{eq:distance-func}. 
Denote by $\tilde M$ the universal cover of $M$ and let $p\from\tilde M\rightarrow M$ be the covering map. Then we have a parametrization of $\tilde X$ by $\tilde M^n$, the Cartesian product of $n$ copies of $\tilde{M}$, as follows.

\begin{proposition}[{\cite[Theorem 2.1]{luo2021deformation}}]
\label{para}
There exists a homeomorphism from the space of geodesic mappings $\tilde X$ to $\tilde M^n$, denoted by $\phi\mapsto\tilde\phi = (\tilde\phi_1,\dots,\tilde\phi_n)$, such that $p(\tilde\phi_i)=\phi(i)$ for any vertex $i\in V$. 
\end{proposition}

We remark that this parametrization is not unique.
For each vertex $i\in V$, we arbitrarily fix a lift of $i$ in $\tilde M$. Then $\tilde \phi$ is determined by the positions of these lifts in $\tilde M$. See \cite{luo2021deformation} for more details on the construction of this homeomorphism.
Hence, $X$ is a $2n$-dimensional open manifold in $\tilde X$. Denote by $\partial X$ the boundary of $X$ in $\tilde X$. Intuitively,  an element of $\partial X$ contains degenerate geodesic triangles but no overlap between the interiors of geodesic triangles.

\subsection{Organization}
This paper is organized as follows. In Section \ref{sec:preliminary}, we summarize prerequisites for proving Theorem \ref{main2}, including a generalized Tutte's embedding theorem established in \cite{luo2021deformation}, the mean value coordinates introduced by Floater in \cite{floater2003mean}, and the construction of an exhaustion of $X$ by compact subsets. 
In Section \ref{sec:main-proof}, we prove Theorem \ref{main2} assuming a key lemma on the existence of a special homotopy in $X$. In Section \ref{sec:key-lemma}, we describe the detailed construction of the homotopy in $X$, with the proofs of some auxiliary propositions deferred to Section \ref{sec:auxiliary}.

\subsection{Acknowledgments.}
The work of Z. R. is supported in part by NSF grant DMS-2603932.

\section{Preparations for the Main Theorem}\label{sec:preliminary}
This section collects three known results, including a generalized Tutte's embedding theorem, the mean value coordinates, and a characterization of subsets of $\mathbb{R}^k$ homeomorphic to $\mathbb{R}^k$.

\subsection{A generalized Tutte's embedding theorem}
Tutte's embedding theorem~\cite{tutte1963draw} is a fundamental result on the construction of straight-line embeddings of a $3$-connected planar graph in the plane, which has profound applications in computational geometry. The connection between Tutte's embedding theorem and spaces of geodesic triangulations of surfaces was established in \cite{luo2022spaces} when the surface is a convex polygon in the plane. This connection was further generalized to negatively curved surfaces
in \cite{luo2021deformation} to prove Theorem \ref{contractible}. We summarize this result briefly below.

Let $(i, j)$ be the directed edge starting from the vertex $i$ and ending at the vertex $j$. Denote by
$
\vec E=\{(i,j):ij\in E\}
$
the set of directed edges of $T$. A positive vector $w\in\mathbb R^{\vec E}_{>0}$ is called a \emph{weight} of $T$. For any weight $w$ and any geodesic mapping $\varphi\in  X$, we call $\varphi$  \textit{$w$-balanced} if for any $i\in V$,
$$
\sum_{j:(i,j)\in \vec E}w_{ij}{v}_{ij}=0.
$$
% $$
% \sum_{j:ij\in E}w_{ij}{v}_{ij}=0.
% $$
Here $v_{ij}=\varphi'_{ij}(0)\in TM$ is the tangent vector of $\varphi_{ij}$ at $\varphi_{ij}(0)$ pointing toward $\varphi_{ij}(1)$.
% Here $ v_{ij}\in T_{q_i}M$ is defined with the exponential map $\exp:TM\to M$ such that $\exp_{q_i}(t{v}_{ij}) = \varphi_{ij}(t)$ for $t\in[0, 1]$, where $q_i$ is a point on $M$.

\begin{theorem}
\label{Tutte}
 For any weight $w$, there exists a unique geodesic triangulation $\varphi\in X(M, T,\psi)$ that is $w$-balanced. Moreover, the induced map $\Phi\from \mathbb R^{\vec E}_{>0}\to X$ given by $\Phi(w)=\varphi$ is continuous.
% $\Phi(w)=\varphi$ is continuous from $\mathbb R^{\vec E}_{>0}$ to $X$.
\end{theorem}

Theorem \ref{Tutte} is a generalization of the embedding theorems by Colin de Verdi{\`e}re~\cite[Theorem 2]{de1991comment}
% (see Theorem 2 in \cite{de1991comment}) 
and Hass--Scott~\cite[Lemma 10.13]{hass2012simplicial}
% (see Lemma 10.12 in \cite{hass2012simplicial}) 
from the case of symmetric weights (i.e., $w_{ij}=w_{ji}$ for any edge $ij$) to non-symmetric weights.

Let 
$$W:=\{w\in \mathbb R^{\vec E}_{>0}:\sum_{j:(i,j)\in \vec E}w_{ij}=1\text{ for any $i\in V$}\}$$
be the space of normalized weights. We then define a map $\alpha\from \mathbb{R}^{\vec E}_{>0}\to W$ by
$$
\alpha(w)_{ij}=\frac{w_{ij}}{\sum_{k:(i,k)\in \vec E}w_{ik}}.
$$
We observe that a geodesic triangulation $\phi$ is $w$-balanced if and only if it is $\alpha(w)$-balanced, i.e.,
\begin{equation}
\label{inclusion}
    \Phi=\Phi\circ\alpha.
\end{equation}

\subsection{Mean value coordinates}
The notion of \emph{mean value coordinates} was introduced by Floater \cite{floater2003mean} as a generalization of barycentric coordinates to produce parametrizations of surfaces in the plane. The same construction works for geodesic triangulations of surfaces. 

Given $\varphi\in X$, the mean value coordinates are defined to be
\begin{equation}
\label{mvc}
w_{ij}=\frac{\tan(\alpha_{ij}/2)+\tan(\beta_{ij}/2)}{l_{ij}},
\end{equation}
where $l_{ij}=| v_{ij}|$ is the geodesic length of $\varphi_{ij}([0,1])$, and $\alpha_{ij}$ and $\beta_{ij}$ are the two inner angles in $\varphi(T^{(1)})$ at $\varphi(i)$ in the two triangles sharing the edge $\varphi_{ij}([0,1])$. This produces a continuous map $\Psi$ from $X$ to $\mathbb R^{\vec E}_{>0}$ by applying the formula (\ref{mvc}) to the pullback of each star of $\varphi(i)$ to the tangent space at $\varphi(i)$. Furthermore, we have $\Phi\circ \Psi=\id_X$ by Floater's mean value theorem (see \cite[Proposition 1]{floater2003mean}). Hence, the mean-value-coordinate map $\Psi$ is a section of the map $\Phi\from \R^{\vec E}_{>0}\to X$.
% from $X$ to $\mathbb R^{\vec E}_{>0}$. 

\subsection{Simple connectivity at infinity}
There are contractible open subsets of $\mathbb{R}^k$ that are not homeomorphic to $\mathbb{R}^k$, such as the Whitehead manifold. The following theorem characterizes when a homeomorphism to $\mathbb{R}^n$ exists. Recall that a topological space $S$ is  \textit{simply connected at infinity} if for any compact subset $K$ of $S$, there is a compact set $D$ in $S$ containing $K$ so that the induced map
$$ \pi _{1}(S-D) \to \pi _{1}(S-K)$$
is trivial.
\begin{theorem}[\cite{stallings1962piecewise}]
\label{simple}
A contractible open subset of $\mathbb{R}^k$ with $k\geq 5$ which is simply connected at infinity is homeomorphic to $\mathbb{R}^k$.
\end{theorem}

\section{Proof of the Main Theorem}\label{sec:main-proof}
In this section, we construct an exhaustion of $X$ by compact subsets and use Theorem \ref{simple} to establish Theorem \ref{main2}. Since any triangulation of $M$ has at least three vertices, $X$ has dimension at least $6$, so Theorem \ref{simple} is applicable.

\subsection{Exhaustion of $X$ by compact sets}
For any $\phi\in X$, denote by $\theta^i_{jk}(\phi)$ the inner angle at vertex $i$ in the geodesic triangle $\phi(\triangle ijk)$, and
$$
\theta_m(\phi)=\min\{\theta^i_{jk}(\phi),\theta^j_{ik}(\phi),\theta^k_{ji}(\phi):\triangle ijk\in F\}.
$$
For any $\epsilon>0$, let $K_\epsilon=\{\phi\in X:\theta_m(\phi)\geq\epsilon\}$. The following proposition shows that the sets $\{K_{1/n}\}_n$ form an exhaustion of $X$ by compact sets.
For any $\phi\in  \tilde X$, we still denote by $l_{ij} = l_{ij}(\phi)$ the length of the geodesic segment $\phi_{ij}([0,1])$ in $M$.

\begin{proposition}
\label{exhaustion}
For any $\epsilon>0$, $K_\epsilon$ is compact.
\end{proposition}
\begin{proof}
We show that any sequence $\phi^{(n)}$ in $K_\epsilon$ has a subsequence converging to some point in $K_\epsilon$. Passing to a subsequence, we may assume that $\{l_{ij}(\phi^{(n)})\}$ converges to some $\bar l_{ij}\in[0,\infty]$ for each $ij\in E$. 

If $\bar l_{ij}=\infty$ for some $ij\in E$, then by the hyperbolic law of cosines and the comparison theorems of geodesic triangles, 
$\theta^i_{jk}(\phi^{(n)})\rightarrow0$ or $\theta^j_{ik}(\phi^{(n)})\rightarrow0$ 
as $n\to \infty$, for any $\triangle ijk\in F$.
This contradicts the fact that $\theta_m(\phi^{(n)}) \geq \epsilon$ for all $n$.
% This is impossible because $\phi^{(n)}\in K_\epsilon$ for any $n$. 
Therefore,  $\bar l_{ij}<\infty$ for any $ij\in E$.
After passing to a further subsequence, by the compactness of $M$ and finiteness of $T$, we may assume that $\phi^{(n)}$ converges to $\bar\phi\in\tilde X$ with $\bar l_{ij}=l_{ij}(\bar\phi)$.

Now we claim that the $\bar l_{ij}$ are all positive.
If not, some $\bar{l}_{ij}$ vanishes. We note that not all $\bar{l}_{ij}$ can be zero, since $\bar\phi$ is not null-homotopic. Thus, there exists $\triangle ijk\in F$ such that $\bar l_{ij}=0$ and $\bar l_{ik}>0$. By the comparison theorems of geodesic triangles and the hyperbolic law of sines, $\theta^k_{ij}(\phi^{(n)})\rightarrow 0$. This again contradicts the fact that $\theta_m(\phi^{(n)}) \geq \epsilon$ for all $n$.

Finally, it suffices to rule out the case in which  there exists $\triangle ijk\in F$ with three positive edge lengths but a degenerate image. 
% Since we have shown that $\bar l_{ij}$ are all positive, 
Thus we may assume that $\bar l_{ij}=\bar l_{ik}+\bar l_{kj}$. This indicates that $\theta^i_{kj}(\phi^{(n)})\rightarrow 0$, which contradicts $\theta_m(\phi^{(n)}) \geq \epsilon$ for all $n$.
\end{proof}

The following proposition presents a useful estimate to control the minimum angle in the triangulation by conditions on the weights.
\begin{proposition}
\label{weightlimit}
For any $\epsilon>0$, there exists $\delta=\delta(M,T,\psi,\epsilon)>0$ 
such that for any weight $w\in W$, if $w_{ij}+w_{ik}\geq1-\delta$ for some $\triangle ijk\in F$, then $\Phi(w)\notin K_\epsilon$.
\end{proposition}
\begin{proof}
We prove the proposition by contradiction.
Suppose there is a sequence of weights $w^{(n)}$ in $W$ and a sequence of triangles $\triangle i_nj_nk_n$ such that
$w_{i_nj_n}^{(n)}+w_{i_nk_n}^{(n)}\rightarrow1$ as $n\to\infty$, and $\phi^{(n)}:=\Phi(w^{(n)})\in K_\epsilon$ for all positive integers $n$.
% $w_{i_nj_n}^{(n)}+w_{i_nk_n}^{(n)}\rightarrow1$ as $n\to\infty$ and $\Phi(w^{(n)})\in K_\epsilon$ for all positive integer $n$. 
Passing to a subsequence, we may assume that 
\begin{enumerate}
    \item the triangles $\triangle i_nj_nk_n$ are the same triangle $\triangle ijk$,
    \item the sequence $\{w^{(n)}\}$ converges to some $\bar w$ in the closure of $W$ with $\bar w_{ij}+\bar w_{ik}=1$, and
    \item the sequence $\{\phi^{(n)}\}$ converges to some $\bar\phi\in K_\epsilon$.
\end{enumerate}
Since $\phi^{(n)}$ is $w^{(n)}$-balanced, we have 
\begin{equation}
    \sum_{t:(i,t)\in \vec E}w^{(n)}_{it}{v}^{(n)}_{it}=0.
\end{equation}
Here $v^{(n)}_{it}\in TM$ is the tangent vector of $\phi^{(n)}_{it}$ at $\phi^{(n)}_{it}(0)$ pointing toward $\phi^{(n)}_{it}(1)$. After taking the limit, we then have
\begin{equation}
    \sum_{t:(i,t)\in \vec E}\bar w_{it}{\bar v}_{it}=0.
\end{equation}
Here $\bar v_{it}\in TM$ is the tangent vector of $\bar{\phi}_{it}$ at $\bar{\phi}_{it}(0)$ pointing toward $\bar{\phi}_{it}(1)$.
By $\bar w_{ij}+\bar w_{ik}=1$ and $\sum_{t:(i,t)\in \vec E}w_{it}=1$, we know that
\begin{equation}
    \bar w_{ij}{\bar v}_{ij} +\bar w_{ik}{\bar v}_{ik}=0.
\end{equation}
Thus $\bar\phi_{ij},\bar\phi_{ik}$ are on the same geodesic, which implies that $\bar\phi(\triangle ijk)$ degenerates.
This contradicts the assumption that $\bar\phi$ is an embedding.
\end{proof}

\subsection{Contraction of loops at infinity }

Since $X$ is known to be a contractible manifold, by Theorem \ref{simple} it suffices to prove the following fact to complete the proof of Theorem \ref{main2}. 

\begin{proposition}
\label{simpleinfinity}
For any $\epsilon>0$, there exists a positive constant $\delta=\delta(M,T,\psi,\epsilon)<\epsilon$ such that any loop $\gamma$ in $X- K_{\delta}$ is contractible in $X-K_{\epsilon}$.
\end{proposition}

The proof of Proposition \ref{simpleinfinity} is based on the following key homotopy lemma.
\begin{lemma}
\label{keylemma}
Fix a vertex $v\in V$.  For any $\epsilon>0$ there exists a positive constant $\delta=\delta(M,T,\psi,\epsilon)<\epsilon$ such that for any loop $\gamma$ in $X-K_\delta$, there exists a continuous map
$$
H\from S^1\times[0,1]\rightarrow \tilde X
$$ 
such that 
\begin{enumerate}
    \item[(1)]  $H(s,0)=\gamma(s)$ and $H(s,1)\in\partial X$, \text{ for any $s\in S^1$},

    \item[(2)] $H(s,t)\in X-K_\epsilon$, \text{ for any $s\in S^1$ and $t\in[0,1)$}, and

    \item[(3)] $(H(s,t))(i)=(H(s,0))(i)=(\gamma(s))(i)$, for any $(s,t)\in S^1\times[0,1]$ and $i\in V-\{v\}$.
\end{enumerate}
\end{lemma}
The proof of Lemma \ref{keylemma} is postponed to Section \ref{sec:key-lemma}. Intuitively, the above homotopy $H$ moves only the vertex $v$ for triangulations in the loop $\gamma$ so that all the triangulations stay in $X - K_\epsilon$ and become degenerate in the end.

We now show that the occurrence of degenerate geodesic triangles in the homotopy $H$ above can be detected by the mean value coordinates of the family of geodesic triangulations in the homotopy as follows. Recall that $\Psi$ is the map given by mean value coordinates and $\alpha$ is the  normalization of weights around each vertex. 

\begin{lemma}
\label{weighthomotopy}
    For any $1>\kappa>0$, there exists $t\in(0,1)$ such that if $w=\alpha\circ\Psi\circ H(s,t)$ for some $s\in S^1$, then there exists $\triangle ijk\in \sstar(v)$ such that 
$$
w_{ij}+w_{ik}\geq1-\kappa.
$$ 
\end{lemma}
\begin{proof}
    Suppose otherwise. Then we can find a sequence $t_n\rightarrow1$ and a sequence $s_n\in S^1$, such that $w^{(n)}=\alpha\circ\Psi\circ H(s_n,t_n)$ satisfies 
$$
w_{ij}^{(n)}+w_{ik}^{(n)}\leq1-\kappa
$$
for any $\triangle ijk\in \sstar(v)$. 
Let $\phi^{(n)}=H(s_n,t_n)$.
Passing to a subsequence, we may assume that $\{s_n\}$ converges to some 
$\bar s\in S^1$. Then $\{\phi^{(n)}\}$ converges to $\bar\phi:=H(\bar s,1)\in\partial X$, and the sequence of edge lengths $\{l_{ij}^{(n)}\}=\{l_{ij}(\phi^{(n)})\}$ converges to $\bar l_{ij}:=l_{ij}(\bar\phi)$, 
where 
$$
\bar l_{ij}= l_{ij}(H(\bar s,1))=l_{ij}(\gamma(\bar s))\in (0,\infty),
$$ 
for any edge $ij$ that is not incident to the vertex $v$. 
For any $s\in S^1$, the only difference between the geodesic mapping in $H(s, t)$ and $\gamma(s)$ is the position of $v$ in $\gamma(s)(\sstar(v))$. Moreover, since the image $\gamma(S^1)$ for a given loop $\gamma$ is compact, the diameter of $\gamma(s)(\sstar(v))$ is uniformly bounded with respect to $s$. Therefore, the diameter of $H(s,t)(\sstar(v))$ is uniformly bounded, and then the $l_{vi}^{(n)}$ are uniformly bounded. Thus, $\bar l_{vi}< \infty$. 
Denote by $\tilde w^{(n)}=\Psi\circ H(s_n,t_n)$
the mean value coordinates of $\phi^{(n)}$ defined by \eqref{mvc} and $ w^{(n)}=\alpha\circ\Psi\circ H(s_n,t_n)$.

\textbf{Case 1:} $\bar l_{vi}>0$ for any $i\in N_v$. Then the inner angles $\bar\theta^i_{jk}:=\theta^i_{jk}(\bar\phi)$ are well-defined for any $\triangle ijk\in F$, and there exists some inner angle $\bar\theta^i_{jk}=\pi$ since $\bar\phi\in\partial X$. 
Since $H(s,0)$ is an embedding for $s\in S^1$, by our construction we know that $\sstar(i)$ never degenerates to a $1$-dimensional complex for any vertex $i$. Thus $\bar\theta^i_{jk}$ is the only inner angle at $i$ in $\bar\phi$ equal to $\pi$. Hence $\tan(\theta^i_{jk}/2)\to\infty$, and then we have $\tilde w_{ij}^{(n)}\rightarrow\infty$, $\tilde w_{ik}^{(n)}\rightarrow\infty$, and for each $m\in N_i-\{j,k\}$, $\{\tilde w_{im}^{(n)}\}$ converges to some real number. Therefore, after normalization, $w_{ij}^{(n)}+w_{ik}^{(n)}\rightarrow1$, which contradicts the assumption that $w_{ij}^{(n)}+w_{ik}^{(n)}\leq1-\kappa$. See the left figure in Figure \ref{figure1}.

\begin{figure}[htbp] 
\begin{center}
\begin{overpic}[width=0.75\linewidth]{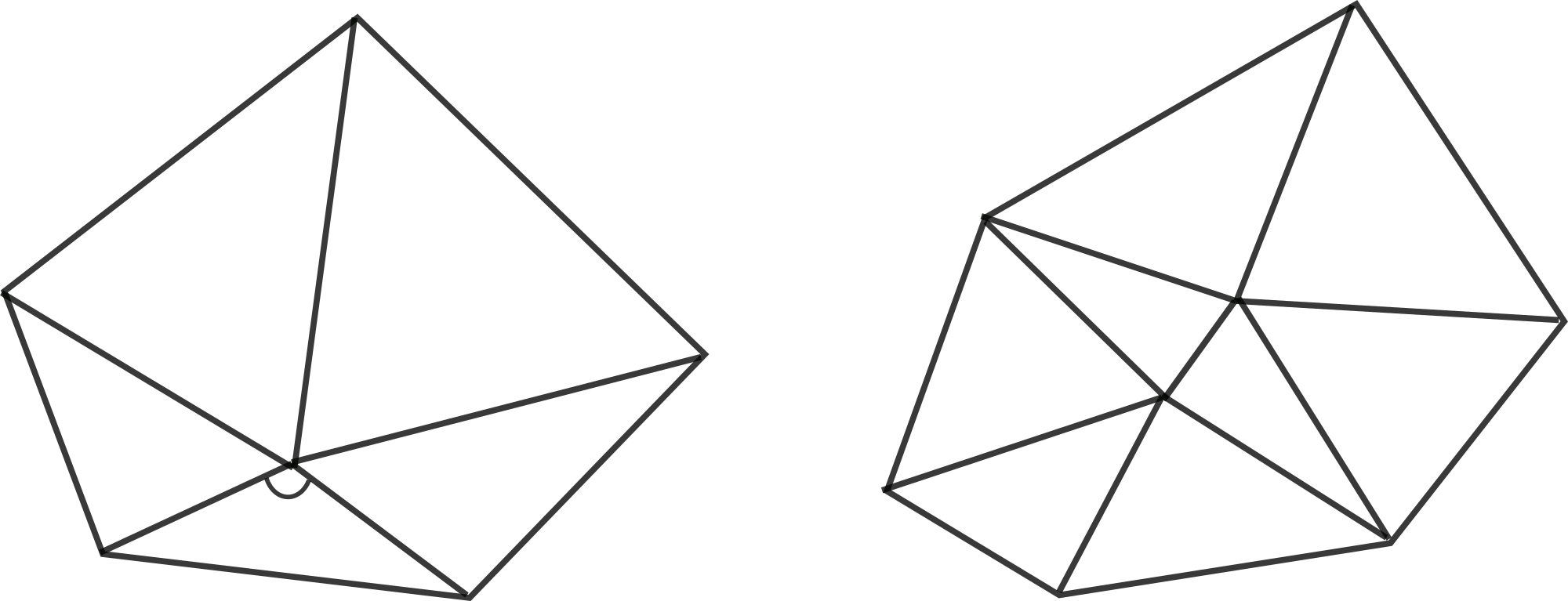}
  \put(20.5,10){$i$}
  \put(5,1){$j$}
  \put(30,-2){$k$}
  \put(16,4){$\approx\pi$}

  \put(60,30){$j=m$}
  \put(80,20){$v$}
  \put(75,10){$i$}
  \put(54,7){$h$}
  \put(92,5){$k$}
\end{overpic}
\caption{Left: Case 1 ($i = v$). Right: Case 2 ($j=m$).}  \label{figure1}
\end{center}
\end{figure}

\textbf{Case 2:} $\bar l_{vi}=0$ for some $i\in N_v$. 
Since $w^{(n)}_{iv}\leq1-\kappa$ for any $n$, by passing to a subsequence, we may assume that there exists some vertex $m\in N_i-\{v\}$ such that for any $n$, 
\begin{equation}
\label{ratio}
\frac{\tilde w_{im}^{(n)}}{\tilde w_{iv}^{(n)}}=\frac{w_{im}^{(n)}}{w_{iv}^{(n)}}\geq\frac{\kappa/(|N_i|-1)}{1-\kappa}.
\end{equation}
Let $j,k$ be the vertices with $\triangle ijv,\triangle ikv\in F$; then by \eqref{mvc}, we have
\begin{equation}
\begin{aligned}
    \tilde w^{(n)}_{iv}&=\frac{1}{l_{iv}^{(n)}}\cdot\left(\tan\frac{\theta^i_{vj}(\phi^{(n)})}{2}+\tan\frac{\theta^i_{vk}(\phi^{(n)})}{2}\right)
\geq\frac{\theta^i_{vj}(\phi^{(n)})+\theta^i_{vk}(\phi^{(n)})}{2l_{iv}^{(n)}}\\
&=\frac{\theta^i_{vj}(H(s_n,t_n))+\theta^i_{vk}(H(s_n,t_n))}{2l_{iv}^{(n)}}
=\frac{\theta^i_{vj}(\gamma(s_n))+\theta^i_{vk}(\gamma(s_n))}{2l_{iv}^{(n)}}\\
&\rightarrow
\frac{\theta^i_{vj}(\gamma(\bar s))+\theta^i_{vk}(\gamma(\bar s))}{0^+}=\infty.
\end{aligned}
\end{equation}
Here, since the vertices $i,j,k$ are fixed in the homotopy $H$ and $H(s,t)\in X$ is an embedding when $t<1$, the sum of angles $$\theta^i_{vj}(H(s,t))+\theta^i_{vk}(H(s,t))=\theta^i_{vj}(\gamma(s))+\theta^i_{vk}(\gamma(s))$$ is  constant with respect to $t$. 
Therefore, Equation (\ref{ratio}) implies 
\begin{equation}
\tilde w^{(n)}_{im}\rightarrow\infty.
\end{equation}

However, if $i,m,v$ do not form a triangle, $\tilde w_{im}^{(n)} = (\Psi(\gamma(s_m)))_{im} \rightarrow (\Psi(\gamma(\bar s)))_{im}<\infty$, since the mean value coordinates remain unchanged if $im$ is not in $\sstar(v)$. Hence,  $m\in \{j, k\}$ and we may assume that $m=j$. See the right picture in Figure \ref{figure1}.

Let $\triangle ijh$ be the triangle in $F$ with $h\neq v$. Notice that by Property (3) of the homotopy in Lemma \ref{keylemma},
$$
\lim_{n\to\infty} l_{ij}^{(n)}=\bar l_{ij}
=l_{ij}((\gamma(\bar s))(i),(\gamma(\bar s))(j))>0,
$$
and
$$
\lim_{n\to\infty} \theta^i_{jh}(\phi^{(n)})=\lim_{n\to\infty} \theta^i_{jh}(\gamma(s_n))=\theta^i_{jh}(\gamma(\bar s))\in(0,\pi).
$$
Then 
$$
\lim_{n\to\infty}\frac{\tilde w_{ij}^{(n)}}{\tilde w_{iv}^{(n)}}
=\lim_{n\to\infty}\frac{1}{\tilde w_{iv}^{(n)}}
\cdot
\frac{1}{l_{ij}^{(n)}}\cdot\left(\tan\frac{\theta^i_{vj}(\phi^{(n)})}{2}+\tan\frac{\theta^i_{jh}(\phi^{(n)})}{2}\right) =\lim_{n\to\infty}\frac{1}{\tilde w_{iv}^{(n)}}
\cdot
\frac{1}{l_{ij}^{(n)}}\cdot\tan\frac{\theta^i_{vj}(\phi^{(n)})}{2}
$$

$$
\leq
\lim_{n\to\infty}
\frac{1}{l_{ij}^{(n)}}\cdot\tan\frac{\theta^i_{vj}(\phi^{(n)})}{2}
\bigg/
\left(
\frac{1}{l_{iv}^{(n)}}\cdot\tan\frac{\theta^i_{vj}(\phi^{(n)})}{2}
\right) =\lim_{n\to\infty}\frac{l_{iv}^{(n)}}{l_{ij}^{(n)}}=0.
$$
This contradicts Equation \eqref{ratio}, which completes the proof. 
\end{proof}

Notice that in Lemma \ref{weighthomotopy}, the choice of $t$ is uniform in $s\in S^1$.
Now we are ready to prove Proposition \ref{simpleinfinity}. Recall that $\Phi$ is the map given by the generalized Tutte's embedding theorem.

\begin{proof}[Proof of Proposition \ref{simpleinfinity}]
 For any $\epsilon>0$, let $\delta$, $\gamma$, and $H$ be as in Lemma \ref{keylemma}. It is straightforward to check that $\alpha\circ\Psi\circ H$ is a continuous map on $S^1\times[0,1)$, and  $\theta_m\circ H(s,t)\rightarrow0$ as $t\rightarrow1$. 
%More specifically, for any $\epsilon_1>0$ there exists $\delta_1>0$ such that if $t\geq1-\delta_1$ then $\theta_m\circ H(s,t)<\epsilon_1$ for any $s\in S^1$. 

By Proposition \ref{weightlimit}, there exists $\kappa>0$ such that for any weight $w$, 
if $w_{ij}+w_{ik}\geq1-\kappa$ for some $\triangle ijk\in F$,
then
$$
\theta_{m}(\Phi(w))<\epsilon.
$$
For this particular $\kappa>0$, Lemma \ref{weighthomotopy} implies that there exists $t_0\in(0,1)$ such that if $w=\alpha\circ\Psi\circ H(s,t_0)$ for some $s\in S^1$, then there exists $\triangle ijk\in \sstar(v)$ such that 
$$
w_{ij}+w_{ik}\geq1-\kappa.
$$ 
By Lemma \ref{keylemma}, $\gamma$ is homotopic to $\beta:=H(\cdot,t_0)$ in $X-K_\epsilon$. Let 
$$
\tilde\beta=\alpha\circ\Psi\circ\beta\from S^1\rightarrow W_0,
$$ 
where 
$$
W_0:=\{w\in W:w_{ij}+w_{ik}\geq1-\kappa\text{ for some $\triangle ijk\in \sstar(v)$}\}.
$$
Since $\Phi\circ\tilde\beta=\beta$, 
it suffices to show that $\tilde\beta$ is null-homotopic in $W_0$ so that $\beta$ is null-homotopic in $X-K_\epsilon$.

Fix a vertex $u\in N_v$. Define $\tilde\pi\from S^1\rightarrow W_0$ 
to be a loop in $W_0$ such that
$$
\tilde\pi_{ij}(s)=
\begin{cases}
\tilde\beta_{ij}(s), & i\neq u,\\[4pt]
1-\kappa, & i=u \text{ and } j=v,\\[4pt]
\dfrac{\kappa}{|N_u|-1}, & i=u \text{ and } j\in N_u\setminus\{v\},
\end{cases}
$$
Now, given any $s\in S^1$, since $\tilde\beta(s)\in W_0$,  $\tilde\beta_{ij}(s)+\tilde\beta_{ik}(s)\geq1-\kappa$ for some $\triangle ijk\in \sstar(v)$. Let $\zeta(s)$ be the line segment in $W$ that connects $\tilde\beta(s)$ and $\tilde\pi(s)$.
\begin{enumerate}
    \item If $u\neq i$, then $\tilde\pi_{ij}(s)=\tilde\beta_{ij}(s)$ and $\tilde\pi_{ik}(s)=\tilde\beta_{ik}(s)$. Hence, for any $w\in\zeta(s)$, $w_{ij}=\tilde\beta_{ij}(s)$ and $w_{ik}=\tilde\beta_{ik}(s)$.
    Thus $w_{ij}+w_{ik}\geq1-\kappa$, and then $w\in W_0$. So we have $\zeta(s)\subset W_0$.
    \item If $u=i$, notice that $v$ is either $j$ or $k$ in the triangle $\triangle ijk$, and we have $$\tilde\pi_{ij}(s)+\tilde\pi_{ik}(s)=\tilde\pi_{uj}(s)+\tilde\pi_{uk}(s)\geq \tilde\pi_{uv}(s)=1-\kappa.$$
    Then for any $w\in\zeta(s)$, $w_{ij}+w_{ik}\geq1-\kappa$, and therefore $\zeta(s)\subset W_0$.
\end{enumerate}
Therefore, we always have $\zeta(s)\subset W_0$, so we know $\tilde\beta$ is linearly homotopic to $\tilde\pi$ in $W_0$. On the other hand, $\tilde \pi$ lies in the following subspace
$$W_u:=\{w\in \mathbb R^{\vec E}_{>0}:\sum_{j:j\sim i}w_{ij}=1\text{ for any $i\in V$, and $w_{ui}=\tilde\pi_{ui}(s)$
for any $i\in N_u$ } \}\subset W_0.$$
Here, we note that $\tilde\pi_{ui}(s)$ is a constant with respect to $s\in S^1$ for any $i\in N_u$.
Hence, $W_u$ is convex and therefore contractible, which implies that $\tilde \pi$ and  hence $\tilde \beta$ is null-homotopic in $W_0$. This completes the proof of Proposition \ref{simpleinfinity}.
\end{proof}

\section{Construction of the Homotopy in Lemma \ref{keylemma}}\label{sec:key-lemma}
We first introduce a few more notions and auxiliary propositions in Section \ref{subsec:notations}, and then prove Lemma \ref{keylemma} in Section \ref{subsec:key-lemma}. The auxiliary propositions are proved in Section \ref{sec:auxiliary}. 

\subsection{Convex Disks and kernels}\label{subsec:notations}

We recall the following definitions. A subset $D$ of $M$ is a \textit{disk} if it is homeomorphic to the closed unit  disk. A disk $D\subset M$ is (strongly) \textit{convex} if for any $x,y\in D$, there exists a unique geodesic arc in $D$ connecting $x$ and $y$. Denote this geodesic by $g=g[D,x,y]\from[0,1]\rightarrow M$, parameterized at constant speed such that $g(0)=x$ and $g(1)=y$. For a subset $D$ of $M$, let the \emph{kernel} $\ker(D)$ of $D$ be the subset of $D$ consisting of all points $x\in D$ such that, for any $y\in D$, there exists a geodesic arc $\gamma$ in $D$ connecting $x$ and $y$. Intuitively, $\ker(D)$ consists of all points from which every other point in $D$ is visible via a geodesic segment. These points are called \textit{eyes} of $D$, and we know that $\ker(D)$ is convex. See \cite{kilicman2014note} for more details.

For the fixed vertex $v$ and any $\phi\in X$, denote by $D(\phi)$ the closed (polygonal) disk in $M$ bounded by the union of the geodesic arcs $\phi_{ij}([0,1])$ with $\triangle ijv\in F$. In other words, $D(\phi)=\phi(\sstar(v))$. Since $\phi$ is originally defined on $T^{(1)}$, we sometimes extend $\phi$ to a homeomorphism from $T$ to $M$, still denoted by $\phi$. 
Since $\phi$ is an embedding, we know that $\phi(v)$ is an eye of $D(\phi)$. Moreover, the points in a small neighborhood of $\phi(v)$ are all eyes of $D(\phi)$.
Therefore, $\ker(D(\phi))$ is a convex disk.

For every $\phi\in X$, $x\in\ker(D(\phi))$, and $y\in D(\phi)$, there exists a unique geodesic arc in $D(\phi)$ connecting $x$ and $y$. We still denote this geodesic by $g=g[D,x,y]\from[0,1]\rightarrow M$, parametrized at constant speed such that $g(0)=x$ and $g(1)=y$. The uniqueness follows from the fact that $D(\phi)\subset M$ is simply connected and equipped with a metric of negative curvature. 

We summarize the necessary properties of convex disks and related maps in the following three propositions. Their proofs are postponed to Section \ref{sec:auxiliary}.

\begin{proposition}\label{prop:projection}
Suppose $D$ is a convex disk in $M$. For any two distinct points $x,y\in D$ with $x\in \inte(D)$, there exists a unique point $z=z(D,x,y)\in\partial D$ such that $y$ is on the geodesic arc $g[D,x,z]$.

\end{proposition}
The map $y\mapsto z(D,x,y)$ is a geodesic radial projection onto $\partial D$ centered at $x$.

\begin{proposition}[\cite{karcher1977riemannian}]
\label{prop:barycenter}
If $D\subset M$ is a convex disk, then there exists a unique minimizer $b(D)\in D$ of the function
$$
E(x)=\int_D(d_D(x,y))^2dA(y),
$$
where $dA$ is the Riemannian area form on $M$. Furthermore, $b(D)\in \inte(D)$.
\end{proposition}
This minimizer is the Riemannian barycenter of mass or the Karcher mean of the domain $D$. See \cite{karcher1977riemannian, karcher2014} for details.

For the fixed vertex $v$, given $\phi\in X$ and $x\in\ker(D(\phi))$, define $\phi_x\from T^{(1)}\rightarrow M$ by
\begin{equation}\label{eq:pertubation-map}
    (\phi_x)_{ij}=\begin{cases} 
      \phi_{ij}, & \text{if } v\notin\{i,j\}, \\
      g[D(\phi),x,\phi(j)], & \text{if } v=i.
\end{cases}
\end{equation}
Notice that $\phi(v)\in \inte(\ker(D(\phi)))$ since $\phi$ is an embedding. Geometrically, $\phi_x$ is the map obtained by perturbing $\phi$ so that $v$ is mapped to $x$ instead of $\phi(v)$. We note that $\phi_x$ is always an element of $\Tilde{X}$, and $\phi_x\in X$ when $x\in\inte(\ker(D(\phi)))$.

Denote by $K = K(M)$ the set of nonempty compact subsets of $M$, which is naturally a metric space under the Hausdorff distance. Let 
$$
K_c=\{D\in K:\text{$D$ is a convex disk}
\},
$$
$$
P=\{(D,x,y)\in K_c\times M\times M:x\in D,y\in D\},
$$
and
$$
P_0=\{(D,x,y)\in P:x\neq y, x\in\inte(D) \}.
$$
We then have 
\begin{proposition} 
\label{convergence}
The following properties hold.

\begin{enumerate}
    \item $D\mapsto\partial D$ is a continuous map from $K_c$ to $K$,
    \item $D\mapsto b(D)$ is a continuous map from $K_c$ to $M$,
    \item $(D,x,y,t)\mapsto g[D,x,y](t)$ is a continuous map from $P\times[0,1]$ to $M$,
    \item $(D,x,y)\mapsto z(D,x,y)$ is a continuous map from $P_0$ to $M$,
    \item $\phi\mapsto\ker(D(\phi))$ is a continuous map from $X$ to $K$, and
    \item $(\phi,x)\mapsto\phi_x$ is a continuous map from $\{(\phi,x)\in X\times M:x\in\ker(D(\phi))\}$ to $\tilde X$, and $\phi_x\in X$ if and only if $x\in \inte(\ker(D(\phi)))$.
\end{enumerate}
\end{proposition}

\subsection{Proof of Lemma \ref{keylemma}}\label{subsec:key-lemma}
\begin{figure}[htbp]  
\begin{center}
\begin{overpic}[width=0.4\linewidth]{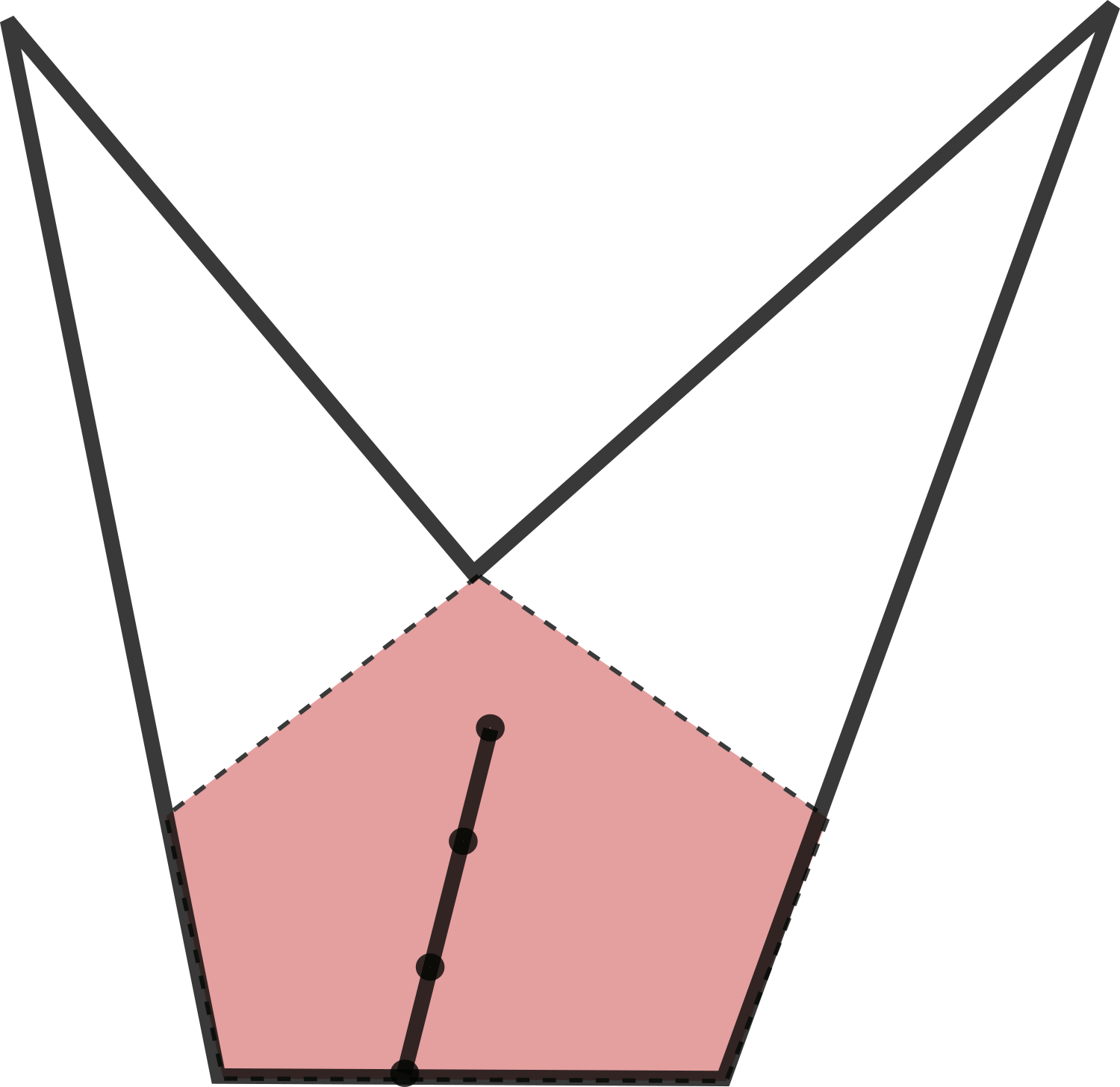}
  \put(45,28){$b(D)$}
  \put(43,20){$\phi(v)$}
  \put(30,-4){$x$}
  \put(41,11){$y$}
  \put(22,5){$D$}
  \put(68,5){$D(\phi)$}

\end{overpic}
\caption{The set $\tilde K_\epsilon$ where the kernel $D$ of $D(\phi)$ is the pink region. }  \label{figure2}
\end{center}
\end{figure}

\begin{proof}[Proof of Lemma \ref{keylemma}]

For simplicity, we let $D := \ker(D(\phi))$, which depends on $\phi$, in this paragraph. Notice that the interior of $D$ is not empty since it contains $\phi(v)$. For any $\epsilon > 0$, consider a subset $\tilde{K}_{\epsilon}$ of $X$ defined as follows: A geodesic triangulation $\phi$ belongs to $ \tilde{K}_{\epsilon}$ if there exist $x \in \partial D$ and $y \in g[D, \phi(v), x]([0, 1])$ such that $\phi(v) \in g[D, b(D), x]([0, 1])$ and $\phi_y \in K_\epsilon$. Formally, it can be written as 
$$
\tilde K_{\epsilon}=\{\phi\in  X:
\text{$\exists x\in\partial D,y\in g[D,\phi(v),x]([0,1])$ s.t. $\phi(v)\in g[D,b(D),x]([0,1])$ and $\phi_y\in K_\epsilon$}
\}.
$$
We note that $K_\epsilon\subset\tilde K_\epsilon$, by taking $y=\phi(v)$ and $x$ be the projection of $\phi(v)$ from $b(D)$ onto $\partial D$. Essentially, the set $\tilde K_{\epsilon}$ consists of geodesic triangulations that can be deformed to $K_\epsilon$ by moving the image of $v$ along a ray starting from $b(D)$. In other words, if $\phi\in X - \tilde K_\epsilon$,  for any $x\in \partial D$, $\phi_y$ does not belong to $K_\epsilon$ for any $y$ in the geodesic segment $g[D, \phi(v), x]$. We then can deform the geodesic triangulation $\phi$ by moving $\phi(v)$ along $g[D, \phi(v), x]$  to $x$ such that $\phi$ stays away from $K_\epsilon$ during this deformation. See Figure \ref{figure2} for an illustration.

We present the following claims and prove them later.

\textbf{Claim 1}: $\tilde K_\epsilon$ is compact for any $\epsilon>0$.

Then for any $\epsilon>0$, we can pick a $\delta>0$ such that $\tilde K_\epsilon\subset K_\delta$.
Let $\gamma$ be a loop in $X-K_\delta$. For simplicity, we let $A :=\sstar(v)$ in the rest of the proof.

\textbf{Claim 2}: The loop $\gamma$ is freely homotopic to $\gamma_0$ in $X-K_\delta$ such that for all $s\in S^1$,
\begin{equation}
\label{pert}
    (\gamma_0(s))(v)\neq b(\ker(\gamma_0(s))(A)).
\end{equation}

Assuming the above two claims and $\gamma=\gamma_0$ without loss of generality, we construct the required homotopy $H(s, t)$ as follows. 
For any $s\in S^1$,
\begin{enumerate}
    \item Let $D_s=\ker((\gamma(s))(A))$ be a convex disk. By Proposition \ref{prop:barycenter},  $b_s=b(D_s)$;
    \item Denote the image of the vertex $v$ by $v_s=(\gamma(s))(v)$. Then by Proposition \ref{prop:projection}, let $z_s=z(D_s,b_s,v_s)$;
    \item With the geodesic segment $y(s,t)=g[D_s,v_s,z_s](t)$ and Equation  \eqref{eq:pertubation-map}, we define the homotopy $H(s,t)=(\gamma(s))_{y(s,t)}$.
\end{enumerate}
Then Property (1) and  (3) in Lemma \ref{keylemma} follow directly from the definition. 
We note that when $t\in[0,1)$, $y(s,t)\in\inte(D_s)$, and hence $H(s,t)=(\gamma(s))_{y(s,t)}\in X$. Moreover, since $\gamma(s) \in X - K_\delta \subset X - \tilde K_\epsilon$, by the definition of $\tilde K_\epsilon$, we know that when $y(s,t)\in\inte(D_s)$, $(\gamma(s))_{y(s,t)}\notin K_\epsilon$. Thus, for any $s\in S^1$ and $t\in[0,1)$, $H(s,t)\in X-K_\epsilon$, and then Property (2) holds.

\textbf{Proof of Claim 1}:
Suppose $\{\phi^{(n)}\}$ is a sequence in $\tilde K_\epsilon$ and let $D_n=\ker(\phi^{(n)}(A))$. Then there exist $x_n\in\partial D_n$, $t_n\in[0,1]$, and $y_n=g[D_n,\phi^{(n)}(v),x_n](t_n)$, such that 
\begin{enumerate}
    \item $\phi^{(n)}(v)= g[D_n,b(D_n),x_n](t_n')$ for some $t_n'\in[0,1]$, and
    \item $\phi^{(n)}_{y_n}\in K_\epsilon$.
\end{enumerate}
Passing to a subsequence, we may assume that $y_n\rightarrow\bar y\in M$, and  $\phi^{(n)}_{y_n}\rightarrow \phi_0\in K_\epsilon$, since $M$ and $K_\epsilon$ are compact. 
We need to prove  $\phi^{(n)}\rightarrow\bar\phi\in \tilde X$ by verifying that the length of each fixed edge of $\phi^{(n)}$ converge to a finite number after passing to a subsequence.
Indeed, note that the edges of $\phi^{(n)}$ and those of $\phi^{(n)}_{y_n}$ differ only at the edges incident to $v$. 
For any edge $ij$ not adjacent to $v$, the length of that edge in $\phi^{(n)}$ is equal to the edge length in $\phi^{(n)}_{y_n}$, which converges to the edge length in $\phi_0$. 
For any edge $vi$ adjacent to $v$, $\phi^{(n)}(vi)$ and $\phi^{(n)}_{y_n}(vi)$ are contained in the same set $\phi^{(n)}(A)$.
Note that the diameters of $\phi^{(n)}(A)$ are uniformly bounded since $\phi^{(n)}_{y_n}\to \phi_0\in K_\epsilon$. 
Hence the length of $\phi^{(n)}(vi)$ is uniformly bounded with respect to $n$. Hence it converges to a finite number after passing to a subsequence.

By Proposition \ref{convergence}, we have 
$$
D_n=\ker(\phi^{(n)}(A))=\ker(\phi^{(n)}_{y_n}(A))\rightarrow\ker(\phi_0(A)):=\bar D,
$$ as $n\to\infty$. Moreover,  $\bar D$ is a convex disk, and $b(D_n)\rightarrow b(\bar D)$.
Passing to a subsequence again, we may assume that $x_n\rightarrow \bar x\in\partial \bar D$, $t_n\rightarrow\bar t\in[0,1]$,
and $t_n'\rightarrow\bar t'\in[0,1]$. Then $\bar y=g[\bar D,\bar\phi(v),\bar x](\bar t)\in D$ and $\bar\phi(v)=g[\bar D, b(\bar D), \bar x](\bar t')$ by Proposition \ref{convergence}. It remains to prove that $\bar\phi\in X$. If not, we have $\bar\phi(v)\in\partial\bar D$. Then $\bar\phi(v)=\bar x$ and $\bar y=\bar x=\bar\phi(v)\in\partial \bar D$. 
This contradicts the fact that $\phi_0\in K_\epsilon$. 

\textbf{Proof of Claim 2}:
% Recall that by Proposition \ref{para}, $\tilde X$ is homeomorphic to $\tilde M^{n}$, and $p\from\tilde M \to M$ is a covering map. 
% We fix a parametrization of $\tilde M$ by $\mathbb R^2$, and then there exists some small $\eta>0$ such that if $a\in\mathbb R^2$ and $\|a\|_2<\eta$, 
% then $\gamma(s)_{p(\widetilde{\gamma(s)}(v)+a)}\in X-K_\delta$ for any $s\in S^1$. Indeed, if not, there exists $s_n\in S^1$ and $a_n\in \mathbb{R}^2$ such that $\|a_n\|_2 \to 0$ but $\gamma(s_n)_{p(\widetilde{\gamma(s_n)}(v)+a_n)}\in K_\delta$. Since $K_\delta$ is compact, by passing to a subsequence, we have $s_n\to \bar s$ and $\gamma(s_n)_{p(\widetilde{\gamma(s_n)}(v)+a_n)}\to  \gamma(\bar s)_{p(\widetilde{\gamma(\bar s)}(v))} = \gamma(\bar s)\in K_\delta$. This is a contradiction. 
Recall that by Proposition \ref{para}, $\tilde X$ is homeomorphic to $\tilde M^{n}$, and $p\from\tilde M \to M$ is a covering map. 
We fix a parametrization of $\tilde M$ by $\mathbb R^2$.
Since $X$ is open in $\tilde{X}$ and $K_\delta$ is compact, $X-K_\delta$ is open in $\tilde{X}$. 
Then by $\gamma(S^1)$ compact, there exists $\eta>0$ such that, for every $s\in S^2$, perturbing the lifted position of the vertex $v$ by any $a\in\R^2$ with $\|a\|_2<\eta$ yields a geodesic triangulation in $X-K_\delta$.
Geometrically speaking, if the local perturbation $a$ of the image $\widetilde{\gamma(s)}(v)$ of the vertex $v$ lifted to $\tilde M$ is small enough, then the geodesic triangulations stay in $X - K_\delta$.

It suffices to construct a continuous function $g\from S^1\rightarrow \mathbb R^2$ such that for any $s\in S^1$, $\|g(s)\|_2<\eta$ and
$$
g(s)\neq(\widetilde{\gamma(s)_{b_s}})(v)-\widetilde{\gamma(s)}(v):=g_0(s).
$$
Then $\gamma(s)_{p(\widetilde{\gamma(s)}(v)+t\cdot g(s))}$
gives a desired homotopy satisfying \eqref{pert}, by the definition of $g(s)$.

Since $g_0$ is continuous, let $g_1\from S^1\rightarrow\mathbb R^2$ be a smooth map such that $\|g_0(s)-g_1(s)\|_2<\eta/2$
for any $s\in S^1$. Since every point in $S^{1}$ is a critical point of $g_1$, Sard's theorem implies that $g_1(S^{1})$ has measure zero in $\mathbb{R}^2$. Then there exists some $a\in\mathbb R^2$ with $\|a\|_2<\eta/2$ such that $g_1(s)\neq a$ for any $s\in S^1$. Therefore,  $g(s)=g_0(s)-g_1(s)+a$ is the required function. 
\end{proof}

\section{Proof of the Auxiliary Propositions}\label{sec:auxiliary}

Since $M$ is a connected closed orientable smooth surface with a smooth Riemannian metric of negative Gaussian curvature, its universal cover $\tilde{M}$ is a $2$-dimensional Hadamard manifold. We will use concepts and properties of Hadamard manifolds including tangent cones, supporting geodesics and  normal vectors of convex disks, and the uniqueness of geodesics between two points.
See \cite{sakai1996riemannian} for more information about Hadamard manifolds.

\subsection{Proof of Proposition \ref{prop:projection}} \label{subsec:proofof4.1}
Since $D$ is simply connected, it suffices to prove the statement for a lift $\tilde{D}$ of $D$ in the universal cover $\tilde{M}$. Let $\tilde{x}$ and $\tilde{y}$ be the lifts of $x$ and $y$ in $\tilde{D}$; then $\tilde{x}\in \inte(\tilde{D})$. Note that the geodesic segment $g[\tilde{D},\tilde{x},\tilde{y}](t)$ can be extended to all $t\geq 0$ in $\tilde M$, and we still denote the geodesic by $g[\tilde{D},\tilde{x},\tilde{y}]$. Consider the following set
\begin{equation}
    S=\{t\geq0:g[\tilde{D},\tilde{x},\tilde{y}](t)\in \tilde{D}\}.
\end{equation}
Since $\tilde{D}$ is compact and convex in $\tilde{M}$, $S$ is a nonempty closed and bounded set. Let $t_0$ be the maximal element of $S$, and $z=g[\tilde{D},\tilde{x},\tilde{y}](t_0)\in\partial\tilde{D}$. It remains to prove the uniqueness.

Suppose $z'\in\partial \tilde{D}$ also satisfies the requirement. Then there is a unique geodesic $l\subset \tilde{M}$ passing through $\tilde{x}$ and $\tilde{y}$. By the assumption, we have $z,z'\in l$. Moreover, there exists $t'_0\in[0,t_0]$ such that $z'=g[\tilde{D},\tilde{x},\tilde{y}](t'_0)$. 

To prove the uniqueness, it suffices to show that for any $t\in[0,t_0)$, $g[\tilde{D},\tilde{x},\tilde{y}](t)\in\inte(\tilde{D})$. Indeed, since $\tilde{x}\in\inte(\tilde{D})$, there is an open neighborhood $N(\tilde{x})$ of $\tilde{x}$ such that $\tilde{x}\in N(\tilde{x})\subset \inte(\tilde{D})$. By the convexity of $\tilde{D}$, for any $w\in N(\tilde{x})$, there is a unique geodesic arc $g[\tilde{D},w,z]\subset\tilde{D}$ that connects $w$ and $z$. Let $U=\cup_{w\in N(\tilde{x})} g[\tilde{D},w,z]\subset\tilde{D}$ be the union of these geodesic arcs.
Since $\tilde{M}$ is a Hadamard manifold, the exponential map $\exp_z\from T_z\tilde{M}\to\tilde{M}$ at $z$ is a diffeomorphism with $\exp_z(\textbf{0})=z$.
It is then sufficient to show that $\exp_z^{-1}(g[\tilde{D},\tilde{x},\tilde{y}](t))$ is an interior point of $\exp_z^{-1}(U)=\cup_{w\in N(\tilde{x})} \exp_z^{-1}(g[\tilde{D},w,z])$ for any $t\in(0, t_0)$.
Since $\exp_z^{-1}(g[\tilde{D},w,z])$ is always a line segment with $\textbf{0}$ as an endpoint, we know that $\exp_z^{-1}(U)$ is the union of line segments that connect $\textbf{0}$ and points in $\exp_z^{-1}(N(\tilde{x}))$, which is an open neighborhood of $\exp_z^{-1}(\tilde{x})$. Then the result follows  from Euclidean geometry in $T_z\tilde M$: the open cone over $\exp_z^{-1}(\tilde{x})$ with the apex $\textbf{0}$ contains  $\exp_z^{-1}(g[\tilde{D},\tilde{x},\tilde{y}](t))$ in the interior of $\exp_z^{-1}(U)$ for $t\in(0, t_0)$.

\subsection{Proof of Proposition \ref{prop:barycenter}}
As in Section \ref{subsec:proofof4.1}, we can lift the convex disk $D$ to $\tilde{D}$ in the universal cover $\tilde{M}$ of $M$, which is a Hadamard manifold, and it suffices to prove the result for $\tilde{D}$. For simplicity, we still denote it by $D$.

It is easy to show that $E(x)$ is continuous and thus has a minimum on $D$. Suppose $x,x'$ are distinct minimizers of $E$, and let $x''=g[D,x,x'](1/2)$. Then for any $y\in D$, by the CAT(0) midpoint inequality for the geodesic triangle $\triangle yxx'$,
$$
(d_D(x,y))^2+(d_D(x',y))^2\geq2(d_D(x'',y))^2+(d_D(x,x'))^2/2>2(d_D(x'',y))^2.
$$
Here, since \(D\) is convex, its intrinsic distance agrees with the ambient distance between the corresponding lifts in \(\tilde M\). Then
$$
2E(x)=E(x)+E(x')=\int_D((d_D(x,y))^2+(d_D(x',y))^2)dA
$$
$$
>\int_D2d_D(x'',y)^2dA=2E(x'').
$$
This contradicts the fact that $E(x)$ is the minimum. Hence, $E(x)$ has a unique minimizer in $D$.

If the minimizer $x$ is on the boundary of $D$, since $D$ is convex, the tangent cone $C_x$ of $D$ at $x$ is a closed convex cone with nonempty interior in $T_x\tilde{M}$ of angle at most $\pi$. Let $v\in T_x\tilde{M}$ be the unit vector in the direction of the bisector of $C_x$. Then $v$ points toward the interior of $D$, and for any $y\in D$, the tangent vector at $x$ pointing toward $y$ makes an angle at most $\pi/2$ with $v$.
Let $\gamma(t)=\exp_x(tv)$, $0\leq t\leq t_0$, where $t_0>0$ is sufficiently small so that $\gamma(t_0)\in\inte(D)$. There exists a neighborhood $V$ of $\gamma(t_0)$ with positive area in $D$ such that for any $y\in V$, the tangent vector at $x$ pointing toward $y$ makes an angle less than $\pi/2$ with $v$. Hence, for $y\in V$,
\begin{equation}
    \left.\frac{d}{dt}\right|_{t=0^+} d(\gamma(t),y)^2 = -2\langle v,\exp_x^{-1}(y)\rangle<0.
\end{equation}
Therefore,
\begin{equation}
    \left.\frac{d}{dt}\right|_{t=0^+} E(\gamma(t)) = \int_D -2\langle v,\exp_x^{-1}(y)\rangle dA<0,
\end{equation}
which contradicts the fact that $x$ is the minimizer.

\subsection{Proof of Proposition \ref{convergence}}

Recall that the Hausdorff distance on the set $K$ is defined as 
$$d_H(X, Y) = \max \{\max_{x\in X} d(x, Y), \max_{y\in Y} d(X, y)\},$$
where $X, Y \in K(M)$ and $d$ is the distance function induced by the Riemannian metric on $M$. Equivalently, 
$$d_H(X, Y)  = \inf\{\epsilon\geq 0 | X\subset Y_\epsilon \text{ and } Y\subset X_\epsilon\},$$
where $S_\epsilon$ is the $\epsilon$-neighborhood of a subset $S$ of $M$. We denote the convergence of sets in the Hausdorff distance by $D_n\to D$. We also use the following characterization of the Hausdorff convergence. Every limit of a convergent sequence \(x_n\in D_n\) lies in $D$, and for every point \(x\in D\), there exists \(x_n\in D_n\) such that \(x_n\to x\).

Whenever $D_n\to D$ occurs in the remainder of the proof, we work inside a fixed compact disk neighborhood $U\subset M$ 
such that $D_n\subset U$ for sufficiently large $n>0$. We lift $D_n$, $D$, and $U$ to the Hadamard universal cover $\tilde M$ and consider $D_n$ and $D$ as convex sets in $\tilde M$ without changing notation. 
All convex sets are then considered as subsets of a fixed compact disk $U$ in the Hadamard universal cover $\tilde M$.

However, due to the topology of $M$, the restriction of the ambient distance function $d_M$ to $D$ can be different from the intrinsic distance function $d_{\tilde M}$ on a lift of $D$ to $\tilde{M}$, where $d_{\tilde M}$ denotes the path metric induced by the Riemannian metric of $\tilde{M}$.
In particular, the restriction of $d_M$ to $D$ may not be the path metric determined by lengths of paths contained in $D$.
Nevertheless, the following lemma shows that we can use the two metrics to define the Hausdorff convergence interchangeably.

\begin{lemma}\label{lem:hausdorff-lift-equivalence}
Let $p\colon\tilde M\to M$ be a Riemannian covering.  Let $U\subset M$ be a compact topological disk, and let $\tilde U$ be a component of $p^{-1}(U)$.  Equip $U$ with the restriction of the Riemannian distance $d_M$ on $M$, and equip $\tilde U$ with the restriction of the Riemannian distance $d_{\tilde M}$ on $\tilde M$.

For every nonempty compact subset $A\subset U$, let $\tilde A=(p|_{\tilde U})^{-1}(A).$ Then, for nonempty compact subsets $A_n,A\subset U$, $d_{H,M}(A_n,A)\to 0$ if and only if $d_{H,\tilde M}(\tilde A_n,\tilde A)\to 0$.
Here the two Hausdorff distances $d_{H,M}$ and $d_{H,\tilde M}$ are defined using $d_M$ and $d_{\tilde M}$, respectively.
\end{lemma}

\begin{proof}
Because $U$ is a simply connected compact disk, $\tilde U$ is compact and $f:=p|_{\tilde U}\colon\tilde U\to U$ and its inverse $f^{-1}$ are uniformly continuous.

Suppose first that $d_{H,M}(A_n,A)\to 0$.  Given $\epsilon>0$, the uniform continuity of $f^{-1}$ provides $\delta>0$ such that  for any $a,b\in U$, if $d_M(a,b)<\delta$, then $d_{\tilde M}(f^{-1}(a),f^{-1}(b))<\epsilon$.
For sufficiently large $n$, we have $A_n\subset (A)^{M}_{\delta}$ and $A\subset (A_n)^{M}_{\delta}$, where $(A_n)^{M}_{\delta}$ is the $\delta$ neighborhood of $A_n$ with respect to the metric $d_M$. It implies that $\tilde A_n\subset (\tilde A)^{\tilde M}_{\epsilon}$ and $\tilde A\subset (\tilde A_n)^{\tilde M}_{\epsilon}$.
Thus $d_{H,\tilde M}(\tilde A_n,\tilde A)<\epsilon$ for all sufficiently large $n$. 

Conversely, suppose that $d_{H,\tilde M}(\tilde A_n,\tilde A)\to0$.  Since the covering map is a local isometry, it does not increase the distance, namely, for any $\tilde x,\tilde y\in\tilde M$, $d_M(p(\tilde x),p(\tilde y)) \leq d_{\tilde M}(\tilde x,\tilde y)$. It then follows directly from the definition of Hausdorff distance that
$d_{H,M}(A_n,A)\leq d_{H,\tilde M}(\tilde A_n,\tilde A),$ and therefore $d_{H,M}(A_n,A)\to0$.
\end{proof}

Therefore, we only need to prove the convergence results in Proposition \ref{convergence} after lifting to the universal cover $\tilde{M}$ with the metric $d_{\tilde M}$. Before proving Proposition \ref{convergence}, we state an elementary fact that will be used repeatedly.  

\begin{lemma}\label{lem:convex-stability}
Let $C_n$ and $C$ be compact convex disks in a Hadamard surface and suppose that $C_n\to C$ in the Hausdorff metric. If $q\in\inte(C)$ and $q_n\to q$, then $q_n\in\inte(C_n)$ for all sufficiently large $n$.
\end{lemma}

\begin{proof}
Suppose otherwise. After passing to a subsequence, $q_n\notin\operatorname{int}(C_n)$ for every $n$.

Let $z_n$ be the metric projection of $q_n$ onto $C_n$ such that $d(z_n,q_n) = d(C_n,q_n)$.  If $q_n\in\partial C_n$, then $z_n=q_n$. Since $q_n\to q\in C$ and $C_n\to C$ in the Hausdorff metric, $d(q_n,C_n)\le d(q_n,q)+d(q,C_n)\to0.$ It follows that
$$ d(z_n,q) \le d(z_n,q_n)+d(q_n,q) =d(q_n,C_n)+d(q_n,q)\to 0. $$
Thus $z_n\to q$.

Choose a supporting geodesic $L_n$ of $C_n$ at $z_n$, and let $\nu_n\in T_{z_n}\widetilde M$ be a unit normal pointing toward the exterior half-plane. Then $C_n$ is contained in the closed inward half-plane $H_n^-$ bounded by $L_n$.
After passing to a subsequence, the compactness of the unit tangent bundle over a compact neighborhood of $q$ indicates that
$(z_n,\nu_n)\to (q,\nu)$ for some unit vector $\nu\in T_q\widetilde M$. The supporting geodesics $L_n$ and their inward half-planes $H_n^-$ converge locally to the geodesic $L$ through $q$ orthogonal to $\nu$ and the corresponding closed half-plane $H^-$, respectively.

We claim that $C\subset H^-$. Indeed, for every $y\in C$, by Hausdorff convergence, there are $y_n\in C_n$ with $y_n\to y$. Since $y_n\in H_n^-$, we have $y\in H^-$ by taking the limit. Hence $C\subset H^-$. However, it means that $q$ lies on a supporting geodesic $\partial H^- = L$ of $C$, and this is impossible for the interior point $q\in\inte (C)$.
\end{proof}

Now we are ready to prove Proposition \ref{convergence}.

\vspace{0.5cm}

\noindent\textbf{Part (1)}:
Take $p_n\in\partial D_n$ and suppose, after passing to a subsequence, that $p_n\to p$.  Therefore $p\in D$ by the Hausdorff convergence.  If $p\in\inte(D)$, Lemma \ref{lem:convex-stability}
% , applied uniformly on a small ball about $p$,  
implies that $p_n\in \inte(D_n)$ for sufficiently large $n$, which is a contradiction.  Thus every limit point of $\partial D_n$ belongs to $\partial D$.

Conversely, fix $p\in\partial D$ and choose $q\in\inte(D)$. Let $\rho$ be the geodesic ray starting at $q$ and passing through $p$. Suppose that $p_n\in D_n$ and $q_n\in D_n$ such that $p_n\to p$ and $q_n \to q$ such that $p_n\neq q_n$ for all $n$. 
By Lemma \ref{lem:convex-stability}, we may assume that $q_n\in \inte(D_n)$.
Let $\rho_n$ be the geodesic ray starting at $q_n$ and passing through $p_n$. Then $\rho_n$ intersects $\partial D_n$ at a unique point $r_n$. Since $D_n$ and $D$ are contained in some compact set $U$, after passing to a subsequence if necessary, we may assume $r_n\to r\in \partial D$ by the previous paragraph. The geodesic ray $\rho_n$, considered as a map from $[0, \infty)$, converges to $\rho$ uniformly on every compact subset of $[0, \infty)$. Hence, $r\in \rho\cap \partial D$. The convexity of $D$ and $q\in\inte(D)$ imply $\rho\cap\partial D=\{p\}$, so $r_n\to r=p$.  This completes the proof of Part (1).

\vspace{0.5cm}

\noindent\textbf{Part (2)}:
Suppose that  $D_n\to D$ in $K_c$, and  use the lifts of $D_n$, $D$, and $U$ in $\tilde M$ as indicated above. For $x\in U$, we let
$$E_n(x)=\int_{D_n}d_{\tilde M}(x,y)^2\,dA(y),\qquad E(x)=\int_Dd_{\tilde M}(x,y)^2\,dA(y).$$
Since the $D_n$ and $D$ are convex disks in $\tilde{M}$,
it is not hard to see that if $x\in U\backslash D$, then there exists $x'\in D$ such that $E(x')< E(x)$.
Together with Proposition \ref{prop:barycenter}, it implies that a minimizer $b(D_n)$ of $E_n$ on $U$ lies in $D_n$ and similarly, the minimizer $b(D)$ of $E$ lies in $D$.

We claim that, by Part~(1), for every $\delta>0$ and all sufficiently large $n$, we have
$$D_n\mathbin\triangle D=(D\backslash D_n)\cup(D_n\backslash D)\subset (\partial D)_\delta.$$ 
Indeed, if $x\in D_n-D$, then $x\in U$ lies outside $D$ and  $d(x, \partial D) = d(x, D)<\delta/2$ for sufficiently large $n$. Similarly, if $y\in D-D_n$, then $d(y, \partial D_n) = d(y, D_n)<\delta/2$ and $d(\partial D_n, \partial D)<\delta/2$ for sufficiently large $n$. 

The boundary of a compact convex disk in a smooth surface is rectifiable and has zero area.  Therefore, $\Area((\partial D)_\delta)\to0$ as $\delta\to0$, and  $\operatorname{Area}(D_n\mathbin\triangle D)\to0$.  It follows that
$$\sup_{x\in U}| E_n(x) - E(x)| \leq \int_{U}d^2_{U}(x, y)|\chi_{D_n}(y) - \chi_D(y)|dA(y) \leq \diam(U)^2 \Area(D\triangle D_n)\to 0.$$

Let $b_n=b(D_n)$. Every subsequence of $\{b_n\}$ has a further subsequence converging to some $\bar b\in \tilde M$. The  uniform convergence of $E_n$ to $E$ yields, for any $x\in U$,
$$ E(\bar b)=\lim E_n(b_n)\leq\lim E_n(x)=E(x).$$
Thus $\bar b$ is the unique minimizer of $E$ by Proposition \ref{prop:barycenter}, namely $b(D)$. Hence every convergent subsequence has the same limit $b(D)$, so $b(D_n)\to b(D)$.

\vspace{0.5cm}

\noindent\textbf{Part (3)}:
In the Hadamard surface $\tilde M$, the logarithm map $\exp_{x}^{-1}(y) = v$ is globally defined and smooth from $\tilde{M}\times\tilde{M}$ to the tangent bundle $T\tilde M$, which gives the initial velocity $v$ of the constant-speed geodesic from $x$ to $y$, parametrized on $[0,1]$. Therefore,  $v$ depends continuously on $x$ and $y$. Then the map sending $(D, x, y, t)$ to $\exp_{x}(t\exp^{-1}_{x}(y)) = g[D,x, y](t)$ is continuous.

\vspace{0.5cm}
\noindent\textbf{Part (4)}:
Let $(D_n,x_n,y_n)\to(D,x,y)$ in $P_0$, and set $z_n=z(D_n,x_n,y_n)$ and $ z=z(D,x,y)$. In the lifts $D_n, D,$ and $U$ as above, consider the geodesic rays
$$\gamma_n(s)=\exp_{x_n}\!\left(s\exp_{x_n}^{-1}(y_n)\right),
 \qquad
 \gamma(s)=\exp_x\!\left(s\exp_x^{-1}(y)\right),
 \qquad s\geq0.$$
There are unique $s_n,s_0\geq1$ such that $z_n=\gamma_n(s_n)$ and $z=\gamma(s_0)$. Moreover, on every bounded subinterval of $[0, +\infty)$, $\gamma_n\to\gamma$ uniformly.

Choose $a,b$ sufficiently close to $s_0$ with $a<s_0<b$ so that $p=\gamma(a)\in\inte(D)$ and $q=\gamma(b)\notin D.$ 
We can assume that the segment $\gamma([a,b])$ lies in a convex normal ball about $z$.  
By Lemma~\ref{lem:convex-stability}, $p_n:=\gamma_n(a)\in\inte(D_n)$ for all sufficiently large $n$.
Since $q\notin D$ and $q_n:=\gamma_n(b)\to q$, Hausdorff convergence implies $q_n\notin D_n$ for all sufficiently large $n$.  Hence, $z_n$ lies on the segment from $p_n$ to $q_n$.  This segment converges to the segment from $p$ to $q$ and can be made arbitrarily small by taking $a,b$ sufficiently close to $s_0$.  It follows that $z_n\to z$.

\vspace{0.5cm}
\noindent\textbf{Part (5)}:
We use the following standard description of the kernel to prove Part (5). 
Fix an orientation on $\tilde{M}$.
For a geodesic polygonal disk $Q \in \tilde M$, 
we orient $Q$ by the induced orientation from $\tilde{M}$, and orient $\partial Q$ by the corresponding boundary orientation.
If $e_i$ is an oriented boundary edge, let $H_i(Q)$ be the closed half-plane on the left of the complete geodesic containing $e_i$. Recall from Section \ref{sec:key-lemma} that for any $\phi\in X$, $D(\phi) = \phi(\sstar(v))$. Let $\tilde D(\phi)$ denote a lift of $D(\phi)$ in $\tilde M$, which is a geodesic polygon containing an eye $\phi(v)$.

\begin{lemma}\label{lem:kernel-halfplanes}
\begin{enumerate}
\item For every geodesic polygonal disk $Q \in \tilde M$,
$$\ker(Q)=\bigcap_i H_i(Q).$$

\item If $\phi_n\to\phi$ in $X$, then
$$\bigcap_i H_i(\tilde D(\phi_n))\to \bigcap_i H_i(\tilde D(\phi)) $$
in the Hausdorff metric, where we choose the lifts such that $\tilde D(\phi_n) \to \tilde D(\phi)$ in $\tilde M$.
\end{enumerate}
\end{lemma}

\begin{proof}
The first claim is the usual half-plane characterization of the kernel of a polygon.  Indeed, a point from which all of $Q$ is visible must lie on the inward side of every boundary edge, since the half-planes are convex. Conversely, if a point lies in all the inward half-planes, the geodesic segment from that point to any point of $Q$ cannot cross a boundary edge from the inside to the outside. Hence, the whole segment is contained in $Q$.

For the second claim, the vertices of every boundary edge vary continuously with $\phi$.  Thus the finitely many bounding geodesics, with their chosen inward sides, vary continuously on $\tilde M$.  Let
$$ K_n=\bigcap_iH_i(\tilde D(\phi_n))\quad \text{and} \quad K=\bigcap_iH_i(\tilde D(\phi)).$$
We can define the supporting function $f_n^i$ for $H_i(D(\phi_n))$ and $f^i$ for $H_i(D(\phi))$ such that for $x\in \tilde M$, $x\in H_i(D(\phi))$ if and only if $f^i(x)\geq 0$. Since $\phi_n\to \phi$ and we choose lifts such that $\tilde D(\phi_n) \to \tilde D(\phi)$ in $\tilde M$, we have $H_i(D(\phi_n))\to H_i(D(\phi))$ and $f_n^i \to f^i$ pointwise and uniformly on compact subsets of $\tilde M$ for all $i$.

If $x_n\in K_n$ and $x_n\to x$, the uniform convergence of $f^i_n\to f^i$ on a compact neighborhood of $x$ implies that $f^i_n(x_n)\to f^i(x)$. Then $f^i(x)\geq 0$ for all $i$ and hence $x\in K$.

Conversely, every interior point $x\in \inte(K)$ satisfies all the strict inequalities  $f^i(x)>0$. This implies that $f_n^i(x)>0$ and hence $x$ belongs to $K_n$ for sufficiently large $n$.  Since $K$ is a convex disk, $\overline{\inte(K)}=K$, which completes the proof. 
\end{proof}

Now we prove Part (5). Suppose $\phi_n\to\phi$ in $X$. 
Lemma~\ref{lem:kernel-halfplanes} then implies $\ker(D(\phi_n))\to\ker(D(\phi))$ in the Hausdorff metric.

\vspace{0.5cm}
\noindent\textbf{Part (6)}:
Under the parametrization $\tilde X\cong\tilde M^{|V|}$ in Proposition~\ref{para}, the map $\phi_x$ is obtained by leaving every lifted vertex except $v$ fixed and replacing the lift of $\phi(v)$ by the compatible lift of $x$.  Hence $(\phi,x)\mapsto\phi_x$ is continuous.  
This also follows edge by edge from Part~(3).

Finally, it is straightforward to see  that $\phi_x$ is an embedding if and only if $x\in\inte(\ker(D(\phi)))$ from the geometric picture of $D(\phi)$.  Indeed, if $x\in\inte (\ker(D(\phi)))$, then every triangle incident to $v$ in $\phi_x$ is nondegenerate and has the same orientation as the corresponding triangle of $\phi$. Thus $\phi_x$ is an embedding. Conversely, if $x\in\partial\ker(D(\phi))$, at least one geodesic triangle in its star becomes degenerate, so $\phi_x$ is not an embedding.

\end{document}